\documentclass[aps, twocolumn, prx, groupedaddress, showpacs, floatfix]{revtex4-2}
\usepackage[dvips]{graphicx}
\usepackage{amsmath}
\usepackage{amsfonts}
\usepackage{amssymb}
\usepackage{amsthm}
\usepackage{xcolor}
\usepackage{times}
\usepackage{comment}
\usepackage{hyperref}
\usepackage{siunitx}

\DeclareMathOperator{\Rep}{Re}

\begin{document}
\renewcommand{\topfraction}{0.98}	        
\renewcommand{\bottomfraction}{0.98}        
\setcounter{topnumber}{3}
\setcounter{bottomnumber}{3}
\setcounter{totalnumber}{4}                 
\setcounter{dbltopnumber}{4}                
\renewcommand{\dbltopfraction}{0.98}        
\renewcommand{\textfraction}{0.05}	        
\renewcommand{\floatpagefraction}{0.5}	    
\renewcommand{\dblfloatpagefraction}{0.5}	
\sloppy 
\newcommand{\beq}{\begin{equation}}
\newcommand{\eeq}{\end{equation}}
\newtheorem{rmk}{Remark}
\newcommand{\divg}{\mbox{\rm{div}}\,}
\newcommand{\Divg}{\mbox{\rm{Div}}\,}
\newcommand{\D}  {\displaystyle}
\newcommand{\DS} {\displaystyle}
\def\sca   #1{\mbox{\rm{#1}}{}}
\def\mat   #1{\mbox{\bf #1}{}}
\def\vec   #1{\mbox{\boldmath $#1$}{}}
\def\scas  #1{\mbox{{\scriptsize{${\rm{#1}}$}}}{}}
\def\scaf  #1{\mbox{{\tiny{${\rm{#1}}$}}}{}}
\def\vecs  #1{\mbox{\boldmath{\scriptsize{$#1$}}}{}}
\def\tens  #1{\mbox{\boldmath{\scriptsize{$#1$}}}{}}
\def\ten   #1{\mbox{\boldmath $#1$}{}}
\renewcommand{\C}{\mathbb{C}}
\newcommand{\R}{\mathbb{R}}
\newcommand{\Z}{\mathbb{Z}}
\renewcommand{\T}{\mathbb{T}}
\newcommand{\abs}[1]{\left|#1\right|}
\renewcommand{\Re}{\Rep}
\newcommand{\Lb}{\mathcal{R}}
\newcommand{\sset}[1]{\left\{#1\right\}}

\newcommand{\cc}{\rho}
\newcommand{\e}{\mathrm{e}}
\newcommand{\al}{a}
\newcommand{\be}{b}

\title{Higher-order phase reduction captures delay-dependent synchronization phenomena in physical oscillator networks}
\author{Babette A.~J.~de~Wolff$^{1,2}$, Sagnik Chakraborty$^3$, Istv\'{a}n Z. Kiss$^3$, Christian Bick$^{1,4,5,6}$, Bob W.~Rink$^1$}
\affiliation{
\it $^1$Amsterdam Center for Dynamics and Computation, Department of Mathematics, Vrije Universiteit Amsterdam, the Netherlands\\
\it $^2$Department of Mathematics, Universit\"at Hamburg, Germany \\
\it $^3$Department of Chemistry, Saint Louis University, Missouri, USA \\
\it $^4$Mathematical Insititute, University of Oxford, Oxford, UK\\
\it $^5$Department of Mathematics, University of Exeter, Exeter, UK\\
\it $^6$Institute for Advanced Study, Technische Universit\"at M\"unchen, Garching, Germany}
\begin{abstract} 
\noindent 
Coupled oscillators with time-delayed network interactions are critical to understand synchronization phenomena in many physical systems.
Phase reductions to finite-dimensional phase oscillator networks yield explicit insights into their dynamics.
However, first-order phase approximations —in which the time delay acts as a phase shift— fail to capture the delay-dependence of synchronization.
We develop a systematic approach to derive phase reductions for delay-coupled oscillators to higher order.
Beyond first order, already a second-order phase reduction captures delay-induced synchronization as demonstrated in coupled Stuart–Landau oscillators and experiments with delay-coupled electrochemical oscillators. Our results establish a general mechanism by which time delays reshape synchronization phenomena and reveal intrinsic limitations of widely used reduced phase models, with implications for a broad range of oscillator networks.

\end{abstract}

\maketitle

Time delays in network interactions---for instance induced by finite transmission speed---shape the synchronization behavior of coupled oscillator networks~\cite{Pazo2016,PhysRevLett.94.088101}. 
These are important, for example, to understand brain dynamics and function~\cite{Deco2009,Petkoski2016}.
Varying delay and coupling strength can induce fully synchronized or desynchronized dynamics~\cite{Kim:2001wo,PhysRevLett.92.114102, Kiss2007} as well as more subtle forms of patterns such as phase clusters, chimera states, and slow switching~\cite{Bick:2017dj, Kiss2007, Kori:2008es, Tinsley:2012ef, 10.1038/ncomms8752}.

As a concrete example in a simple network, consider the synchronization dynamics of two delay-coupled oscillators.
Experiments with delay-coupled electrochemical oscillators exhibit either in-phase~\cite{Kiss:2005us} or anti-phase~\cite{Nagao:2016cm} synchronization at weak coupling, with the selected state determined by the effective sign of the delayed interaction.
However, as illustrated in Fig.~1, increasing the coupling strength and delay leads to bistability between in-phase and anti-phase synchronized states. 
This behavior is not specific to the experimental system: Analogous bistability was observed in prototypical nonlinear oscillator models such as the Brusselator~\cite{Prigogine1968} (see Appendix~\ref{app:Supp1}). 

Importantly, in these observations the coupling strength remains sufficiently weak so that amplitude effects---such as amplitude death or oscillator quenching~\cite{KOSESKA2013173}---are absent. These observations thus fall within the validity range of a phase description. 
However, they cannot be captured by traditional first-order phase reduction alone, a widely used dimension-reduction technique to understand such emergent collective phenomena in coupled oscillator networks~\cite{Nakao2015, Ashwin2015, Pietras2019}.
The bistability and its dependence on the time delay is due to the higher-order contributions to the phase dynamics, specifically second harmonics in the interaction function.
As we show below, these can be derived analytically to provide a description of both the model and the experimental observations.

\begin{figure}
    \centering
    \includegraphics[width=\linewidth]{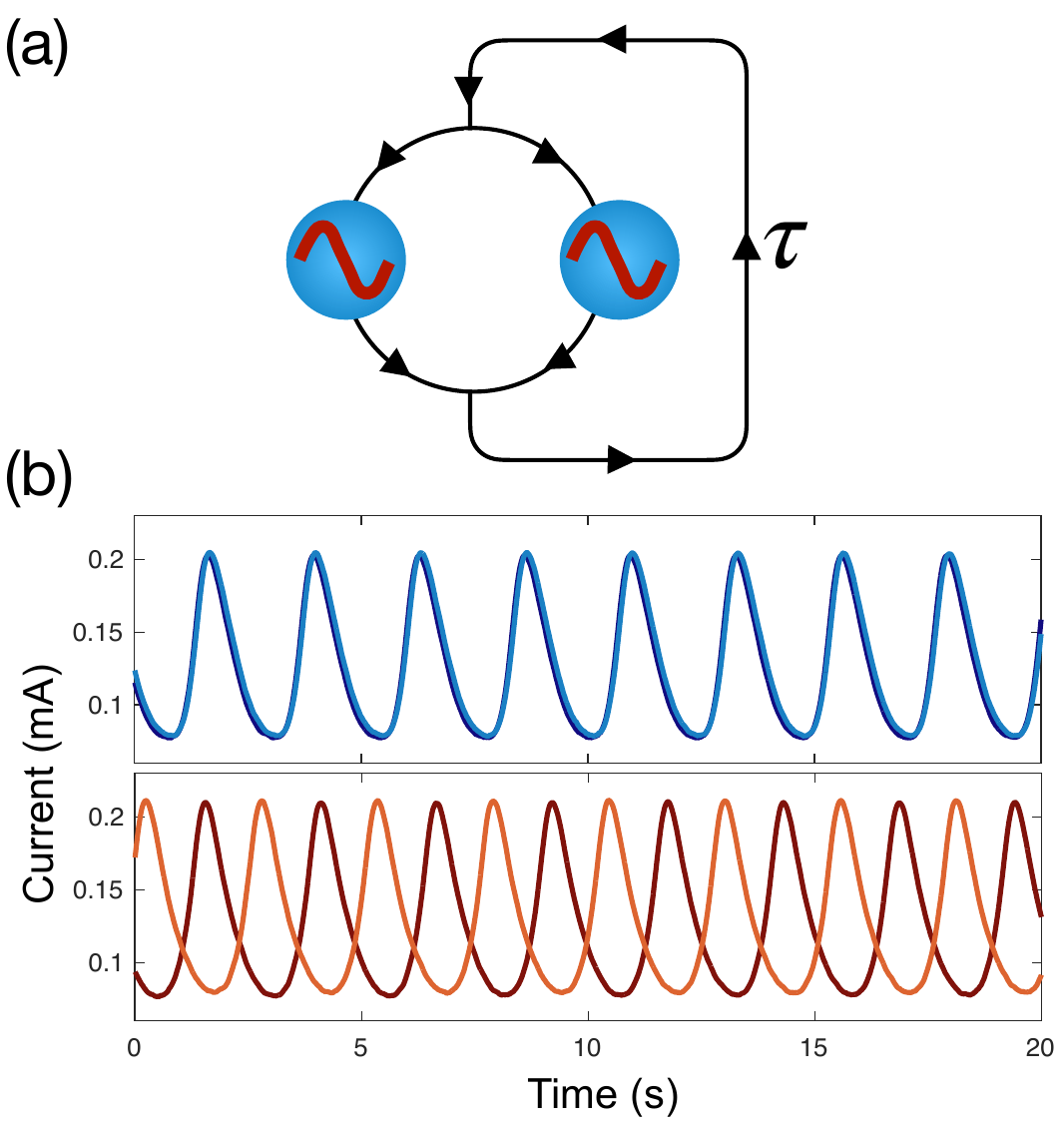}
    \caption{
    Time series of current oscillations of two electrochemical oscillators coupled with a delay.
    (a) Schematic of two oscillators  subject to time-delayed coupling.  
    (b) Current time series of oscillatory nickel electrodissolution with two electrodes coupled with time delay. Experiments with two different initial conditions yield either in-phase (top) or anti-phase (bottom) synchronization. ($K=-125$ V/A$^{-1}$, $\tau=1.85 s$. See Appendix~\ref{app:Supp2} for further experimental details.) Such bistability phenomenon cannot be captured by sinusoidally coupled Kuramoto oscillators with a phase shift, which is a commonly used phase model.}
    \label{fig:timeseries}
\end{figure}

To account for time delays, the standard approach in phase reduction since the seminal work by Kuramoto and Sakaguchi~\cite{Kuramoto,Sakaguchi1986}
is to approximate them by phase lags: 
If $\phi_j(t)\in \mathbb{T} :=\R/2\pi\Z$ is the phase of oscillator~$j$ at time~$t$ that evolves at intrinsic frequency~$\omega$ in isolation then $\phi_j(t-\tau)\approx\phi_j(t)-\omega\tau$.
For coupling strength~$\varepsilon$, this approximation is generally valid as long as $\varepsilon\tau$~is small~\cite{Smirnov2023}.
This limits the applicability of the approximation to real-world physical systems:
In experimental settings, strong coupling may be required to suppress relatively large oscillator heterogeneities or to ensure rapid transient behavior. 
In other cases, the intrinsic oscillator dynamics may be very fast (for example, in lasers~\cite{PhysRevA.45.R4225}), necessitating the use of comparatively large delays~\cite{Yanchuk_2017}.
Moreover, the phase-reduced equations one obtains at first order can be inherently limited:
For example, computing phase reductions for Stuart--Landau oscillators yields interactions with a single harmonic which cannot support  the bistability that is observed in the unreduced system.

In this paper, we develop a systematic approach to higher-order phase reduction for time-delayed systems, and show that higher-order terms clarify how time delay shapes synchronization beyond the phase-shift approximation.
From a mathematical perspective, one can expect that the (infinite-dimensional) dynamics of $n$~delay-coupled oscillators is effectively  described by a finite-dimensional ordinary differential equation on~$\mathbb{T}^n$~\cite{Fenichel1972}.
By contrast, phase-amplitude approximations~\cite{Kotani2020}, second-order Kuramoto dynamics~\cite{Smirnov2025}, or phase dynamics with time delay~\cite{Izhikevich1998} have additional degrees of freedom.
Recent work~\cite{delaypreprint2026} combining a reformulation of a time-delayed system with a parameterization method~\cite{VonderGracht2023a} has made it possible to systematically compute phase reductions to arbitrary order in the presence of coupling delays.
As it computes the phase reduction as a parameterized invariant manifold, it provides information beyond linear stability of individual synchronized states.
It allows us to draw conclusions about global aspects of the dynamics (e.g., basins of attraction), and gives explicit expressions for the amplitudes of the oscillatory dynamics.

The higher-order phase reduction that we derive predicts synchrony and bistability between different synchronized states observed in delay-coupled oscillators.
First, it gives an accurate description of the dynamics of Stuart--Landau oscillators at second order.
Going beyond the phase-shift approximation, the resulting finite-dimensional system correctly predicts  dynamical properties, including the bistability observed in Fig.~1(b), as coupling strength and time delay are varied. 
Second, the  phase dynamics also give a qualitative prediction of the dynamics observed in real-world experiments, including parameter regions of observed bistability and how they change as the time delay increases; cf.~Fig.~1.

\newcommand{\tr}{\mathsf{T}}
\newcommand{\bphi}{\boldsymbol{\phi}}

\section
{Higher-order phase reduction for delay-coupled oscillators}%
To predict synchrony we compute higher-order phase reductions for network connections subject to delay~\cite{delaypreprint2026}.
More specifically, we consider a general model of $n$~nonlinear oscillators coupled with transmission delay. 
In isolation, the state $x_j\in\mathbb{R}^{m_j}$ of oscillator~$j$ evolves according to $\dot x_j(t)=F_j(x_j(t))$ which admits a $T_j$-periodic limit cycle~$x_j^*(t)$ of frequency $\omega_j = \frac{2\pi}{T_j}$.
In other words, the (asymptotic) state of each oscillator is determined by its \emph{phase}~$\phi_j \in \mathbb{T}$ through~$x_j^*(\phi_j / \omega_j)$ on the $j$th limit cycle.
For coupling strength~$\varepsilon\geq 0$, the oscillator $j\in\{1, \dotsc, n\}$ evolves according to 
\begin{align}\label{networkequation}
 \dot x_j(t) = F_j(x_j(t)) + \varepsilon \sum_{k=1}^n G_{j,k}\left(x_j(t), x_k(t-\tau_{j,k})\right),
\end{align}
where the function~$G_{j,k}$ specifies the coupling from oscillator~$k$ to~$j$ and~$\tau_{j,k}\geq 0$ the transmission delay.
We seek to compute the ordinary differential equations that determine the the evolution of the 
$\bphi = (\phi_1, \dotsc, \phi_n)^\tr \in \mathbb{T}^n$
as an expansion in the coupling strength~$\varepsilon \geq 0$ given by
 \begin{align}\label{ansatz2}
 \dot {\bphi} = {\vec f}({\bphi}, \varepsilon) = \vec{\omega} + \varepsilon {\vec f}^{(1)}(\bphi) + \varepsilon^2 {\vec f}^{(2)}(\bphi) +\dotsb,
\end{align}
where~$\vec\omega = (\omega_1, \dotsc, \omega_n)$ are the uncoupled (zeroth order) dynamics and~${\vec f}^{(\ell)}$ the corrections of order~$\ell$.

The first step in our approach is to replace Eq.~\eqref{networkequation}---a DDE with time delays---by an ordinary differential equation coupled with a transport equation. 
Specifically, we consider the co-evolution of oscillator histories~$X_j(s,t)$ and states $x_j(t) = X_j(0, t)$ given by%
\begin{subequations}\label{networkequationcoupled}
\begin{align}
   &\dot x_j(t) = F_j(x_j(t)) + \varepsilon \sum_{k=1}^n G_{j,k}\left(x_j(t), X_k(-\tau_{j,k}, t)\right)\\
      &\partial_t X_j(s,t) = \partial_s X_j(s,t)
   \label{eq:Transport}.
\end{align}
\end{subequations}
Note that problems~\eqref{networkequation} and~\eqref{networkequationcoupled} are equivalent: 
Solutions to the transport equation~\eqref{eq:Transport} with boundary condition $X_j(0,t)=x_j(t)$ are given by $X_j(s,t) = x_j(t+s)$.

In the second step, we match the dynamics of the transformed system~\eqref{networkequationcoupled} to the phase dynamics~\eqref{ansatz2} to be computed.
To this end, we seek functions~$\vec e$, $\vec E$  that relate solutions~$\bphi(t)$ of the phase dynamics~\eqref{ansatz2} to solutions $\vec X = (X_1, \dotsc, X_n)$, $\vec x = (x_1, \dotsc, x_n)$ of~\eqref{networkequationcoupled} 
through
\begin{align}  \label{ansatz}
 \vec x(t) & = \vec e(\bphi(t), \varepsilon), &
    \vec X(s,t) & = \vec E(s,\bphi(t), \varepsilon) 
\end{align}
with $\vec e(\bphi, \varepsilon) = \vec E(0, \bphi, \varepsilon)$. 
For fixed $\varepsilon \geq 0$ and  $\bphi\in \mathbb{T}^n$, the map $s\mapsto \vec E(s,\bphi, \varepsilon)$ describes a ``history'' of the dynamics, that is, a continuous function of the delay parameter $s\in (-\infty, 0]$.
Substitution of~\eqref{ansatz2} and~\eqref{ansatz} into~\eqref{networkequationcoupled}
gives
\begin{subequations}\label{conjeqn}
\begin{align} 
    \partial_{\bphi}  e_j(\bphi, \varepsilon) \cdot {\vec f}(\bphi, \varepsilon)   & =  { F}_j(e_j({\bphi}, \varepsilon)) \nonumber  \\ 
    & \hspace{-2cm} + \, \varepsilon \sum_k G_{j,k}(e_j(\bphi, \varepsilon), E_k(-\tau_{j,k}, \bphi, \varepsilon)),  
    \\ 
    \partial_{\bphi} E_j(s, \bphi, \varepsilon) \cdot {\vec f}(\bphi, \varepsilon)   & = \partial_s E_j(s, \bphi, \varepsilon).
\end{align}
\end{subequations}
Fulfilling these \emph{conjugacy equations} links ${\vec e}$, ${\vec E}$ and the phase dynamics ${\vec f}$ to the dynamics of~\eqref{networkequationcoupled}---and thus~\eqref{networkequation}.

The third step is to expand~\eqref{conjeqn} in powers of~$\varepsilon$ to obtain iterative equations to determine ${\vec e}$, ${\vec E}$, and ${\vec f}$. 
Write
\begin{align*}  \label{ansatzexp}
 \vec e(\bphi(t), \varepsilon) &= \vec e^{(0)}(\bphi(t))+\varepsilon\vec e^{(1)}(\bphi(t))+\dotsb\\
 \vec E(s,\bphi(t), \varepsilon) &= \vec E^{(0)}(s,\bphi(t))+\varepsilon\vec E^{(1)}(s,\bphi(t))+\dotsb
\end{align*}
and substitute with~\eqref{ansatz2} into~\eqref{conjeqn}.
Collecting zeroth-order ($\varepsilon^0$) terms in~\eqref{conjeqn} of order gives
\begin{subequations}\label{eq:Oone}
\begin{align}
\partial_{\bphi} e_j^{(0)}(\bphi) \cdot \vec \omega &= {F_j}(e^{(0)}_j(\bphi))\\
\partial_{\bphi} E_j^{(0)}(s, \bphi)\cdot \vec \omega &= \partial_s E_j^{(0)}(s, \bphi).
\end{align}
\end{subequations}
For $\ell\geq 1$, collecting $\ell$th-order ($\varepsilon^\ell$) terms in~\eqref{conjeqn} yields
\begin{subequations}\label{Oeps}
\begin{align}
\nonumber\partial_{\bphi} e^{(0)}_j(\bphi) \cdot \vec f^{(\ell)}(\bphi)  +  
\partial_{\bphi} e^{(\ell)}_j(\bphi)  \cdot \vec  \omega\hspace{.25cm}&   \\ 
-  DF_j(e^{(0)}_j(\bphi)) \cdot {e_j}^{(\ell)}(\bphi)   &= h_j^{(\ell)}(\bphi) \label{eq:OepsA}\\
 \partial_{\bphi} E_j^{(\ell)}(s, \bphi) \cdot {\vec \omega} - \partial_s E_j^{(\ell)}(s, \bphi) \omega\hspace{.25cm}
 &=  H_j^{(\ell)}(s, \bphi)\label{eq:OepsB}
 \end{align}
 \end{subequations}
with the $(m_j\times m_j)$-Jacobian~$DF_j$ of~$F_j$ and inhomogeneities $h_j^{(\ell)}, H_j^{(\ell)}$ that are functions of ${\vec e^{{(l)}} }$, ${\vec E}^{{(l)}}$, ${\vec f^{{(l)}}}$ with $0\leq l < \ell$.
For $\ell=1$ these are given explicitly by $h_j^{(1)}(\bphi) = \sum_k G_{j,k}\big(e_j^{(0)}(\bphi), E_k^{(0)}(-\tau_{j,k}, \bphi)\big)$ and 
 $H_j^{(1)}(s, \bphi) = - \partial_{\bphi}E_j^{(0)}(s, \bphi) \cdot \vec f^{(1)}(\bphi)$.
Thus, \eqref{Oeps} are inductive equations that determine the phase dynamics to arbitrary order~$\ell$.

The fourth step is now to solve these equations order-by-order.
The zeroth-order equations~\eqref{eq:Oone} have solution $e_j^{(0)}(\bphi) = x_j^*(\phi_j / \omega_j)$,
$E_j^{(0)}(s,\bphi) = e_j^{(0)} (\bphi+\vec\omega s  ) = x_j^*(\phi_j / \omega_j + s)$, i.e., the periodic orbits of the uncoupled system.
For the first-order ($\ell=1$) equations, \eqref{eq:OepsA} is a system of inhomogeneous linear equations for~$\vec e^{(1)}$ and~$\vec f^{(1)}$ that depend on $\vec e^{(0)}$, $\vec E^{(0)}$.
These equations can be solved in exactly the same way as for the problems without time delay, for example using Floquet theory~\cite{VonderGracht2023a}.
Because the inhomogeneous right hand sides~$h_j^{(1)}(\bphi)$ depend on $E_k^{(0)}(-\tau_{j,k}, \bphi) = e_k^{(0)}(\bphi-\vec\omega\tau_{j,k})$, the 
delays~$\tau_{j,k}$ appear in~$\vec e^{(1)}$ and~$\vec f^{(1)}$ as phase lags with respect to the unperturbed periodic orbits as expected.
The first-order equation~\eqref{eq:OepsB} is an inhomogeneous linear transport equations for $\vec E^{(1)}$. 
It can be solved using the method of characteristics, which---with boundary condition $\vec E^{(1)}(0,\bphi) = \vec e^{(1)}(\bphi)$ and with right hand side $H_j^{(1)}(s, \bphi) = - \partial_{\bphi}E_j^{(0)}(s, \bphi) \cdot \vec f^{(1)}(\bphi)$ ---gives 
\begin{equation} \label{characteristics}
\begin{aligned} 
E_j^{(1)}(s,\bphi)  &= e_j^{(1)}(\bphi + \vec \omega s)  \\ & \qquad + \partial_{\phi_j} E_j^{(0)}(\bphi, s) \cdot \int_0^s f_j^{(1)}(\bphi + (s- \sigma) \omega) d \sigma
\end{aligned}
\end{equation}
Equations~\eqref{Oeps} for order $\ell > 1$ have the same structure: They consist of inhomogeneous linear equations for~$\vec e^{(\ell)}$ and~$\vec f^{(\ell)}$ that can be solved with Floquet theory and a system of inhomogeneous transport equations for~$\vec E^{(\ell)}$ that can be solved by the method of characteristics. 
Our method thus allows to analytically compute phase reductions of~\eqref{networkequation} to arbitrary order in the coupling strength~$\varepsilon$. 

\section
{Delay-coupled Stuart--Landau oscillators}%
As a proof of principle, we apply our method to compute the second-order phase reductions for two identical delay-coupled Stuart--Landau~(SL) oscillators with delayed mean-field coupling,
\begin{equation} \label{eq:sl_coupled}
\begin{aligned}
\dot{z}_1 &= (\al + i \be) z_1 - \bar z_1 z_1^2 + \varepsilon \e^{i \cc} \left(z_1(t-\tau) + z_2(t-\tau) \right) \\
\dot{z}_2 &= (\al + i \be) z_2 - \bar z_2 z_2^2 + \varepsilon \e^{i \cc} \left(z_2(t-\tau) + z_1(t-\tau) \right)
\end{aligned}
\end{equation} 
with $z_j\in \mathbb{C}$ and $i:=\sqrt{-1}$. 
We assume that $\al > 0$ and $\be\neq 0$, so that in the absence of coupling, the SL oscillators possess identical limit cycles $z_j^*(t)=\sqrt{\al} \e^{i\be t}$ of frequency $\be = \omega_j =: \omega$. 
As a result we have $e_j^{(0)}(\bphi) = \sqrt{\al} \e^{i\phi_j}$, $E_j^{(0)}(s, \bphi) = \sqrt{\al} \e^{i(\phi_j+\omega s)}$ and $\vec\omega = (\omega,\omega)^\tr$ to order~$\varepsilon^0$.

\subsection{Computing phase reduced dynamics}

To obtain the first-order approximation, the symmetry of the problem allows to avoid the use of Floquet theory, and instead make a \emph{phase-difference ansatz}
\begin{align*}
\vec e^{(1)}(\bphi) &= (e^{i\phi_1}A(\phi_2-\phi_1),e^{i\phi_2}A(\phi_1-\phi_2))^\tr, \nonumber \\ \nonumber 
\vec f^{(1)}(\bphi) & = (B(\phi_2-\phi_1), B(\phi_1-\phi_2))^\tr,
\end{align*}
with real-valued functions~$A, B$ on $\mathbb{T}$.
This transforms~\eqref{Oeps} into the single complex equation 
\begin{align}\label{ABeqn}  
i \sqrt{a} B(\theta) + 2 a A(\theta) = \sqrt{a} (e^{i(\cc - \omega \tau)} + e^{i(\theta - \omega \tau + \cc)})
\end{align}
Taking the imaginary part of this equation gives the expression for ~$B(\theta)$, from where we find that 
\[
\vec f^{(1)}(\bphi) = \begin{pmatrix}
        \sin(\phi_2 - \phi_1 - \omega \tau + \cc) +\sin(\cc - \omega \tau) \\
        \sin(\phi_1 - \phi_2 - \omega \tau + \cc) +\sin(\cc - \omega \tau)
    \end{pmatrix}.
\]
Thus, up to first order, the dynamics for the \emph{phase difference coordinate} $\psi : = \phi_1 - \phi_2$ is given up to first order by 
\begin{equation} \label{psi_ord1}
    \dot{\psi} = -2 \varepsilon \cos(\cc- \omega \tau) \sin(\psi). 
\end{equation}
Taking the real part of \eqref{ABeqn} gives the expression for ~$B(\theta)$ and thus for $\vec e^{(1)}(\bphi)$; next we use~\eqref{characteristics} to calculate that~$\vec E^{(1)}$ is given by 
\begin{align*}
    & \vec E^{(1)}(s, \bphi) = \\
    &\frac{1}{2 \sqrt{\al}} \!\begin{pmatrix}
    \e^{i \phi_1} \left[\cos(\phi_2 - \phi_1 - \omega \tau + \cc) + \cos(\cc - \omega \tau) \right] \\
   \e^{i \phi_2} \left[\cos(\phi_1 - \phi_2 - \omega \tau + \cc) + \cos(\cc - \omega \tau) \right]
\end{pmatrix}
 \\  
 & +  s   \begin{pmatrix}
      i \sqrt{\al} \e^{i(\phi_1 + \omega s)}  \left[\sin(\phi_2 - \phi_1 - \omega \tau + \cc) + \sin(\cc-\omega \tau) \right] \\
     i \sqrt{\al} \e^{i(\phi_2 + \omega s)} \left[\sin(\phi_1 - \phi_2 - \omega \tau + \cc) + \sin(\cc- \omega \tau) \right]
 \end{pmatrix}.
\end{align*}
This gives a nontrivial dependence on the history parameter~$s$ through the second term in the formula.

To compute the second-order corrections~$\vec f^{(2)}$ to the phase dynamics,
note that~$\vec E^{(1)}(-\tau,\bphi)$ enters into the inhomogeneous right hand side of the iterative equations.
We can therefore expect these to depend on the time delay~$\tau$ in a nontrivial fashion.
We solve the second order equation using computer algebra (see below for the code) and find we find that the dynamics of the phase difference coordinate $\psi = \phi_1 - \phi_2$ up to second order is given by 
\begin{equation} \label{psi_ord2}
    \begin{aligned}
        \dot{\psi} = - 2 &\left (\varepsilon \cos(\cc - \omega \tau) - \varepsilon^2 \tau \cos(2(\cc - \omega \tau)) \right) \sin(\psi) \\
    & +  \frac{\varepsilon^2}{2a} \left[ -1 - 2 a \tau + \cos(2(\cc - \omega \tau)) \right] \sin(2 \psi).
    \end{aligned}
\end{equation}

This specifies the phase dynamics~\eqref{ansatz2} and their delay dependence up to second order; we omit explicit expressions of~$\vec e^{(2)}$ and $\vec E^{(2)}$ as they are not relevant for the phase dynamics.

\renewcommand{\alph}{\alpha}

\subsection{Delay-dependent synchronization}%

The phase equations~\eqref{psi_ord2} now allow to elucidate how the synchronization behavior of the network depends on the time delay~$\tau$. In the phase difference coordinates $\psi = \phi_1 - \phi_2$, \emph{in-phase synchrony} corresponds to the equilibrium $\psi=0$ and \emph{anti-phase synchrony} to the equilibrium $\psi=\pi$. Thus, the synchronization dynamics of the delay-coupled oscillator network depends on the stability properties of these equilibria.

\begin{figure*}[t]
    \centering
    \includegraphics[width=\linewidth]{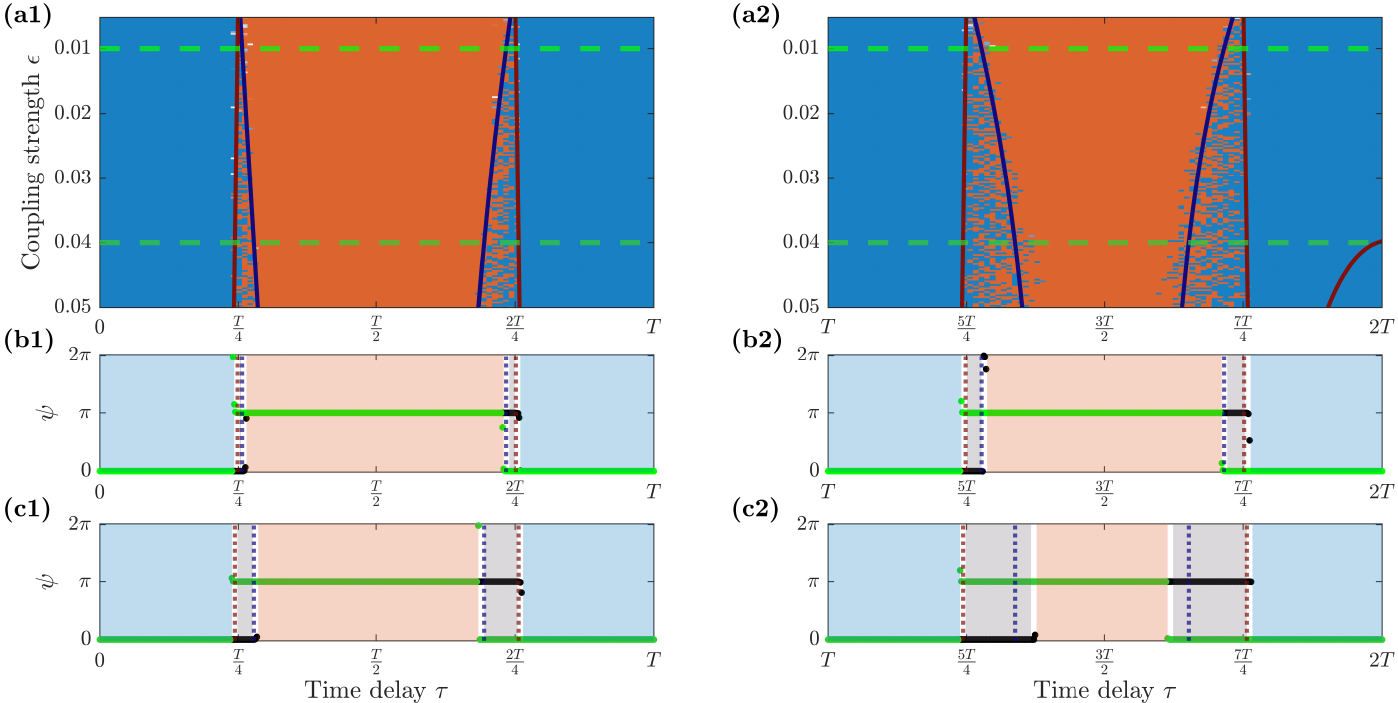}
    \caption{The one-dimensional second-order phase dynamics~\eqref{psi_ord2} predict the dynamics of two delay-coupled Stuart--Landau oscillators~\eqref{eq:sl_coupled}.     
    In all figures, the Stuart--Landau parameters are $\alpha = 1, \beta = 1$ and coupling constant $\rho$ is chosen as $\rho = 0$. 
    The left panels (numbered~1) are for values of the delay up to one period~$T$ of the uncoupled oscillator, the right panels (numbered~2) show results for delays between~$T$ and~$2T$.
    The fact that left and right panels are distinct shows that a phase-shift approximation is not sufficient to capture the synchronization dynamics.
    A spurious bifurcation curve in Panel~(a2) indicates that the second-order phase reduction breaks down as~$\varepsilon \tau$ increases.
    Panels~(a) show that the bifurcations predicted from the second order phase dynamics~\eqref{psi_ord2}---bifurcation of in-phase synchrony as a dark blue line; bifurcation of anti-phase synchrony as a dark red line---predict the bistability observed in the delay-coupled system~\eqref{eq:sl_coupled}:
    The coloring reflects the asymptotic value of $\hat\psi = \arg(z_1\bar z_2)$ ($\hat\psi = 0$ in blue; $\hat\psi = \pi$ in orange).
    The light green line corresponds to the parameter value for the quasi-adiabatic parameter variation in Panels~(b) for small coupling $\varepsilon = 0.01$; and the dark green line corresponds to the parameter value for the quasi-adiabatic parameter variation in Panels ~(c) for larger coupling strength $\varepsilon = 0.04$. 
    In Panels~(b), values of $\hat\psi$ in the forward scan correspond to black dots, the backward scan to light green dots.
    The dashed lines indicate the on/off-set of bistability as predicted by the phase dynamics; the shading indicates mono/bistability of in-phase and anti-phase synchrony.
    Similarly, Panels~(c) show for larger coupling strength $\varepsilon = 0.04$ the forward scan (black dots) and backward scan (dark green dots), showing that larger coupling leads to larger regions of bistability.}
    \label{fig:figure2}
\end{figure*}

To first order the time delay enters the reduction as a phase lag. Specifically, the dynamics for the phase difference truncated at the first order \eqref{psi_ord1} is governed by a single harmonic. Since in \eqref{psi_ord1} the time delay enters as a phase lag, the first-order phase dynamics predicts that the synchronization behavior is fully determined by the effective phase lag~$\alph = \cc- \omega \tau$: 
The equilibrium $\psi=0$ is linearly stable for $\alph\in(-\frac{\pi}{2},\frac{\pi}{2})$, the equilibrium $\psi=\pi$ is linearly stable for $\alph\in(\frac{\pi}{2},\frac{\pi}{2})$, and there is an exchange of stability in a degenerate bifurcation at $\alph\in\{\frac{\pi}{2},\frac{3\pi}{2}\}$.
At the same time, solving the time-delayed differential equations~\eqref{eq:sl_coupled} numerically shows that for reasonable coupling strength the first-order approximation breaks down for time delays well below a single period: The coloring in Fig.~\ref{fig:figure2}(a) has a nontrivial dependence on~$\tau$.

The second-order phase equations \eqref{psi_ord2} do capture the delay-dependence of the synchronization behavior. Note that the second-order terms in \eqref{psi_ord2} contain Fourier coefficients that depend explicitly on the time delay~$\tau$, meaning that  $\tau$ does not enter the reduction as a  phase lag.
 In the second-order phase reduction, linear stability of in-phase synchrony $\psi=0$ and anti-phase synchrony $\psi=\pi$ is determined by
\begin{align}
   \notag 
\lambda^{(0)} & = -\cos(\cc - \omega \tau) - 2 \varepsilon \left( \frac{1}{2a} + \tau \right) \sin^2(\cc - \omega \tau)\, , \\ \notag
\lambda^{(\pi)} & =  \cos(\cc - \omega \tau) - 2 \varepsilon \left( \tau + \left( \frac{1}{2a} - \tau \right) \sin^2(\cc - \omega \tau)  \right)\, .
\end{align}
The bifurcation curves $\lambda^{(0)}=0$ and $\lambda^{(\pi)}=0$ correspond to the dark blue and dark red lines in Fig.~\ref{fig:figure2}(a).
These provide not only a more accurate approximation for the numerical simulations (coloring) but also qualitatively capture the delay-dependence of the synchronization transition.

Naturally, the second-order phase dynamics also capture delay-dependent multistability in the network.
Indeed, simulations of the delay-coupled SL oscillators~\eqref{eq:sl_coupled} with random initial conditions show that there is a large delay-dependent range of parameters for which the system exhibits bistability between in-phase and anti-phase synchrony; this corresponds to the speckled area in Fig.~\ref{fig:figure2}(a).
The region of bistability is analytically predicted to be bounded by the curves $\lambda^{(0)}=0$ and $\lambda^{(\pi)}=0$.  
These curves---determined from the one-dimensional second-order phase reduction and shown as blue/red lines in Fig.~\ref{fig:figure2}---indeed provide approximate bounds for the area of bistability for fixed coupling strength of $\varepsilon=0.1$ up to delays on the order of one period of the oscillation.

\section{Delay-induced Synchronization Pattern with Two Electrochemical Oscillators}

\begin{figure*}[t]
    \centering
    \includegraphics[width=1\linewidth]{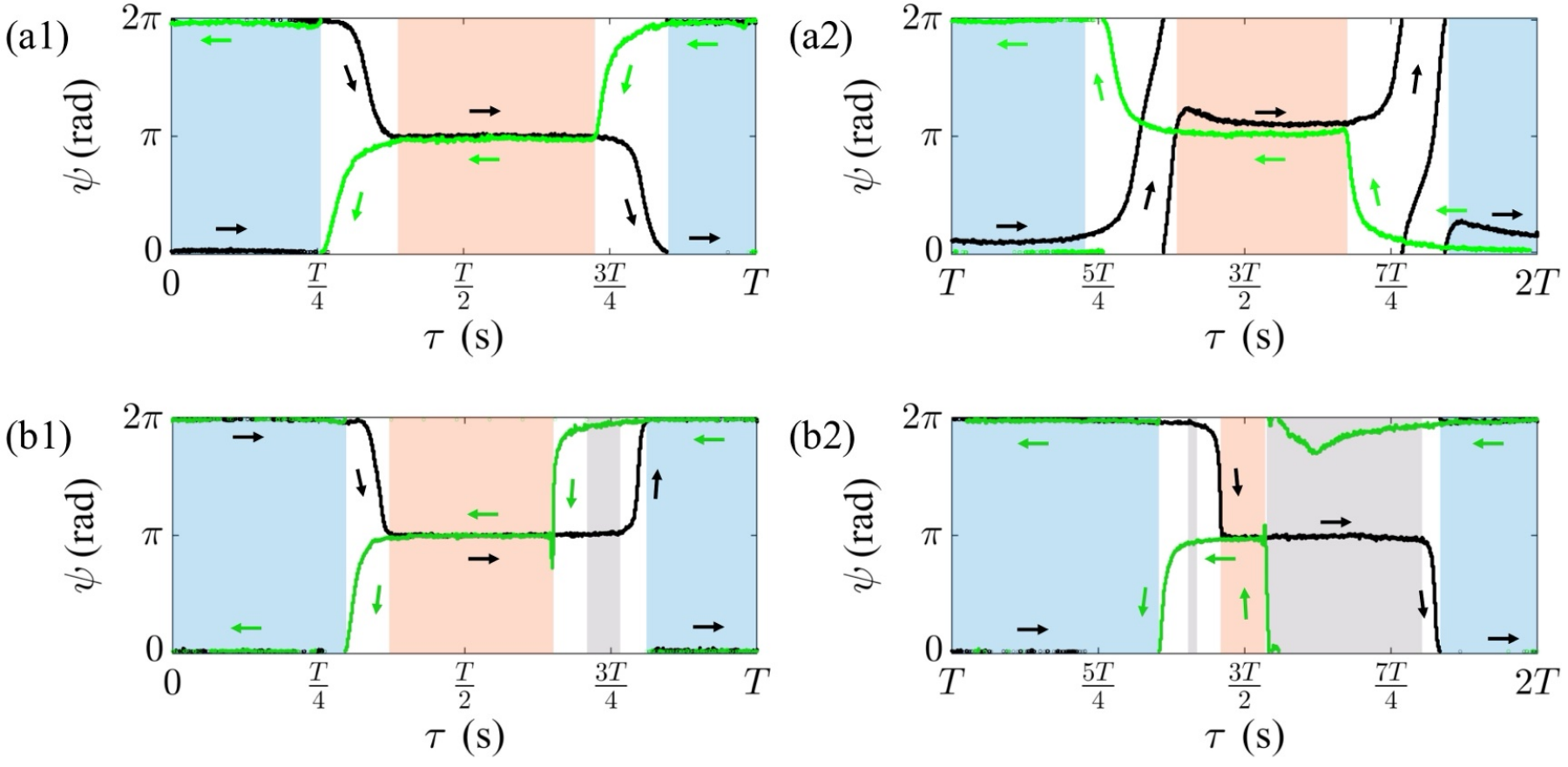}
    \caption{The phase difference~$\psi$ extracted from experiments with two coupled electrochemical oscillators reveals delay-dependent bistability in line with the theoretical predictions.
    (a)~Weak coupling case ($K = -50$ VA$^{-1}$) showing forward (black) and backward (lime) scans for $0 \le \tau \le T$ (a1) and $T \le \tau \le 2T$ (a2). (b) Strong coupling ($K = -125$ VA$^{-1}$) showing forward (black) and backward (lime green) scans when $0 \le \tau \le T$ (b1) and $T \le \tau \le 2T$ (b2). The arrows indicate the direction of changing $\tau$ and the blue, red, and gray shading show unique in-phase, unique anti-phase, or bistable regions, respectively.} 
    \label{fig:1}
\end{figure*}

The qualitative predictions of bistability and delay-dependence of the bifurcations also capture the dynamics of two electrochemical oscillators coupled with time-delayed linear feedback.
The setup (Fig.~\ref{fig:Supp2}) consists of two Ni wires as working electrodes (WE), connected to a potentiostat. The electrodes were immersed in a 3.0 M H\textsubscript{2}SO\textsubscript{4} solution as an electrolyte. 
At a constant circuit potential, with an added external resistance, the electrochemical dissolution of nickel exhibits periodic current~($i_k$) oscillations as a result of the varying chemical composition of the metal-electrolyte interface~\cite{Haim1992ModelingPA}.
The delayed feedback-induced synchronization phenomenon was studied close to Hopf bifurcation at a constant offset circuit potential, $V_{0} = 1080\textrm{ mV}$ where the natural (uncoupled) frequencies of the oscillators were $\omega_1 \approx 405\textrm{ mHz}$ and $\omega_2\approx408\textrm{ mHz}$ , yielding  a mean frequency $\omega_\text{mean}\approx406\textrm{ mHz}$  and mean period $T_\text{mean} = 2.46\textrm{ s}$. 

Coupling is introduced through a real-time controller that implements a global external feedback according to
\begin{equation}
\begin{aligned}
V(t) &= V_0 + K(I(t-\tau) - o)\\
\end{aligned}
\end{equation}
where~$V(t)$ and~$V_0(t)$ are the applied and offset circuit potentials, respectively, $K$~is the coupling strength, $I(t)=i_1(t)+i_2(t)$ is the total current, $\tau$ is a time delay, and $o$~is an offset correction (set to the time averaged total current). 
Without coupling ($K = 0 \textrm{ VA}^{-1}$), the two oscillations are independent and thus the phase difference drifts linearly with time (see Fig.~\ref{fig:Supp3}(a)).

We examined the system dynamics under 
both weak and strong coupling regimes while slowly varying the feedback delay~$\tau$ over two successive delay intervals, $0 \le \tau \le T$ and $T \le \tau \le 2T$---this corresponds to the quasi-adiabatic continuation in Fig.~\ref{fig:figure2}(b,c). 
Each scan included increasing~$\tau$ from $0\textrm{s}$ (or~$T$) to~$T$ (or~$2T$), and then decreasing from~$T$ (or~$2T$) to $0\textrm{s}$ (or~$T$). 
The system was initialized with an in-phase synchronized state, which is a phase locked state of the system with weak coupling at $\tau = 0$, $T$, and $2T$. One such in-phase synchronized state is shown in Fig.~\ref{fig:Supp3}(b). 

The scans were performed very slowly so that each experiment recorded about 1000 cycles. 
The delay-scanning rate (approximately $1\textrm{s}$ delay change over $1000\textrm{s}$ of recording) was chosen to be sufficiently slow so that the observed dynamics represent a quasi-stationary behavior at each value of~$\tau$. 
These experiments clearly capture the occurrence of phase locked oscillatory states at different delays, while also resolving the relatively quick transitions between the states. 

\subsection{Dynamics for weak coupling}
At relatively weak coupling ($K = -50 \textrm{ VA}^{-1}$), in the first delay interval  ($0 \le \tau \le T$) the results are shown in the top panel of Fig.~\ref{fig:1}(a). 
In the forward scan (black), the system  remained in-phase synchronized for small delays. 
As~$\tau$ was increased, a transition occurred to anti-phase synchronization ($\psi \approx \pi$) near $\tau \approx 0.38 T$, before returning to in-phase synchrony near $\tau \approx 0.83T$. A typical anti-phase synchronized state is shown in Fig.~\ref{fig:Supp3}(c). In the backward scan (light green), the oscillators show qualitatively similar behavior.
We thus see that, outside the transition regions, only one stable synchronized state exists for each~$\tau$ in this regime. The blue and the red shading in the figure indicate 
delay region where in- or anti-phase behavior can be clearly identified for both the forward and the backward scans. For example, anti-phase synchronized behavior was observed for $\Delta\tau=0.81$ s ($0.33T$, red shading) during scans in both directions. 
These experiments can be interpreted with delay acting as a phase lag in a Kuramoto oscillator: For small and large delays ($0 \le \tau < T/2$ and $3T/4 < \tau \le T$) in-phase, for intermediate delays ($T/2 < \tau < 3T/4$) anti-phase synchronization could be expected. 
Small deviations do occur from the theoretical predictions, which could be attributed to heterogeneities in oscillator properties, non-stationary behavior (due to slow scan, i.e., transition regions), and slow drift of natural frequencies of the oscillators during the scans. 
 
The bottom panel in Fig.~\ref{fig:1}(a) shows the results for the scan for the next delay interval ($T \le \tau \le 2T$).  
The synchronization regions exhibit similar qualitative features as those seen for the scan  between~$0\textrm{ s}$ and~$T$. 
Two notable differences include a relatively slower transition between the in- and anti-phase synchronized regions and small difference between the phase locked values in the forward and backward scans. Consequently, the delay window with unique anti-phase  synchronization (red shading) shrinks slightly to $\Delta\tau=0.70\textrm{ s}$ ($0.29T$).

\subsection{Dynamics for strong coupling}
Fig.~\ref{fig:1}(b) shows the behavior of the system when coupling is stronger ($K = -125 \textrm{ VA}^{-1}$). While with increasing~$\tau$ there are in-phase through anti-phase to in-phase transitions just like with weak coupling,  the synchronization landscape changes qualitatively. For $0 \le \tau \le T$, (left panel in Fig.~\ref{fig:1}(b)), over an extended range ($0.64T \le \tau \le 0.80T$, grey shading), two distinct phase-locked branches coexist: One near $\psi \approx 0$ (backward scan, dark green) and one near $\psi \approx \pi$ (forward scan, black), exhibiting bistability. Compared to the weak coupling case, the transitions between synchronization states are sharper and broaden into finite delay intervals where both synchronized states are stable. In this case, the anti-phase synchronized window (red) was $\Delta\tau=0.69\textrm{ s}$ ($0.28T$).

The right panel in Fig.~\ref{fig:1}(b) shows the behavior in the next delay interval ($T \le \tau \le 2T$) in the strong coupling regime. Bistable domains emerge, but with noticeable shifts in the location and width of the delay windows. The anti-phase synchronized window (red) shrinks from $\Delta\tau=0.69$ s ($0.28T$) when $0 \le \tau \le T$, to $\Delta\tau=0.20\textrm{ s}$ ($0.08T$) when $T \le \tau \le 2T$. In particular, the bistable regions (gray speckled shading) are significantly wider than in the weak coupling case and even the strong coupling case in the first delay interval, indicating that stronger coupling amplifies the influence of delayed feedback and enhances higher-order delay-induced phase interactions. 

\subsection{Comparison of weak and strong coupling}

We thus observe that while weak coupling yields phase-lag–dominated synchronization with a single phase-locked state at each delay (almost no bistability), the stronger coupling case produces extended, delay-dependent bistability between in-phase and anti-phase synchronization, consistent with higher-order delay-induced phase dynamics, also depicted by the shrinking anti-phase-only regions. 
Notably, in the weak–coupling regime for $0 \le \tau \le T$, the unique anti–phase window spans a relatively large delay interval ($0.33T$), which is comparable to the value expected from the phase–lag approximation ($0.50T$). In contrast, under strong coupling for $T \le \tau \le 2T$, this anti–phase region shrinks to a  narrow window of about $0.08T$. This behavior indicates that small higher–order, delay–induced phase corrections  become dominant as the coupling strength and delay increases.

\section{Discussion}
Higher-order phase reduction of delay-coupled oscillators yields predictions of synchronization phenomena that are in qualitative agreement with experimental observations. 
In particular, the higher-order phase reduction correctly predicts the delay-dependence of bistability regions, which cannot be captured by a phase-lag approximation. 
In the first-order phase-lag approximation, the stability of in-phase and anti-phase synchronization depends only on the effective phase $\alpha = \rho - \omega \tau$ and gives periodic (in $\tau$) stability intervals of each state. 
This behavior is approached for very low coupling strengths, but the experimentally observed narrowing and delay-dependent deformation of the window where only the anti-phase state is stable, reveals deviations for the pure phase-lag picture.
By contrast, higher-order phase reduction theory predicts that second-order corrections introduce additional harmonics and explicit $\tau$-dependent terms in the phase dynamics, which generate delay-dependent bistability regions and shift the synchronization boundaries. 
This is in qualitative agreement with the deformation of stability windows and the emergence of extended bistable regions with increasing coupling and increasing time-delay as measured in the electrochemical oscillator system. 
Note that the prediction agrees despite the real-world oscillators not being exactly identical.
More generally, our approach also allows to directly incorporate small differences between oscillators (beyond leading order) through the self coupling~$G_{k,k}$ in Eq.~\eqref{networkequation}; cf.~\cite{Baibolatov2026}.

Applying our phase reduction approach to larger networks is straightforward and will shed light on how the time delay affects phase interactions to arbitrary order. 
This is particularly relevant when first-order phase reductions are intrinsically limited: 
For example, globally and linearly coupled Stuart--Landau oscillator networks yield a phase reduction with integrable dynamics~\cite{Watanabe1994,Bick2018c} leading to degenerate bifurcations and lack of bistability.
Thus, second-order phase reductions are necessary to uncover the generic dynamical features of the unreduced system through phase reduction.
Note that even if the coupling between nonlinear oscillators is additive (as in~\eqref{networkequation}), the phase interactions beyond first order will have nonpairwise terms that describe indirect phase interactions~\cite{Bick2023c}.
Thus, phase reduction can link phase-lag parameters in higher-order interaction networks~\cite{Bick2021,Battiston2020}---that have so far been chosen ad-hoc---to physically meaningful quantities like time delays.
We will discuss this in more detail in future work. 
Moreover, the approach developed here not only yields the phase dynamics~\eqref{ansatz2} to arbitrary order but through the functions~$\vec e, \vec E$ also an explicit expressions how the dynamics are embedded into the phase space of the full system. This information can be used to make global statements about connecting orbits between different synchrony states or to explain observed amplitude variations~\cite{Ocampo-Espindola2024}.

To predict synchronization phenomena for more complex physical or chemical oscillators, for example relaxation oscillators, one needs more complex models than Stuart--Landau oscillators analyzed here. 
For more complex oscillator models, an analytical computation of the higher-order phase reduced model is at best cumbersome and in some cases unfeasible; therefore, a computational tool that quickly computes the higher-order phase reduced dynamics is needed. 
A numerical implementation of the algorithm from~\cite{delaypreprint2026} and the subsequent analysis of synchronization phenomena for a neural relaxation oscillator is work in progress~\cite{wedgwood202x}. 

\appendix

\section{Simulations with the Brusselator Model}
\label{app:Supp1}

To demonstrate the generality of the findings,  we studied delay-induced synchronization in two coupled Brusselator oscillators, a prototypical nonlinear chemical oscillator model. The equations governing the coupled system are:
\begin{equation}
\begin{aligned}
       \dot{x}_j &= A + x_j^2 y_j - (B_j+1)x_j\\
       &\qquad+ K[x_1(t-\tau) + x_2(t-\tau) - 2o]\\ 
   \dot{y}_j &= B_j x_j - x_j^2 y_j 
\end{aligned}
\end{equation}
where $x_j$ and $y_j$ are the state variables of oscillators $j=1,2$. A small frequency mismatch was introduced through heterogeneous bifurcation parameters $B_1 = 2.09$ and $B_2 = 2.11$ ($A=1$). The natural frequencies of the oscillators are 1.583 x 10$^{-1}$~Hz and 1.581 x 10$^{-1}$~Hz, with  mean  period $T_\text{mean}$ = 6.322~s. The coupling was introduced through the delayed $x$ variable, where $K$ is the coupling strength, and $\tau$ is the delay, and $o=1$ is the offset around which the oscillations take place in the $x$~variable.  

\begin{figure}
    \centering
    \includegraphics[width=\linewidth]{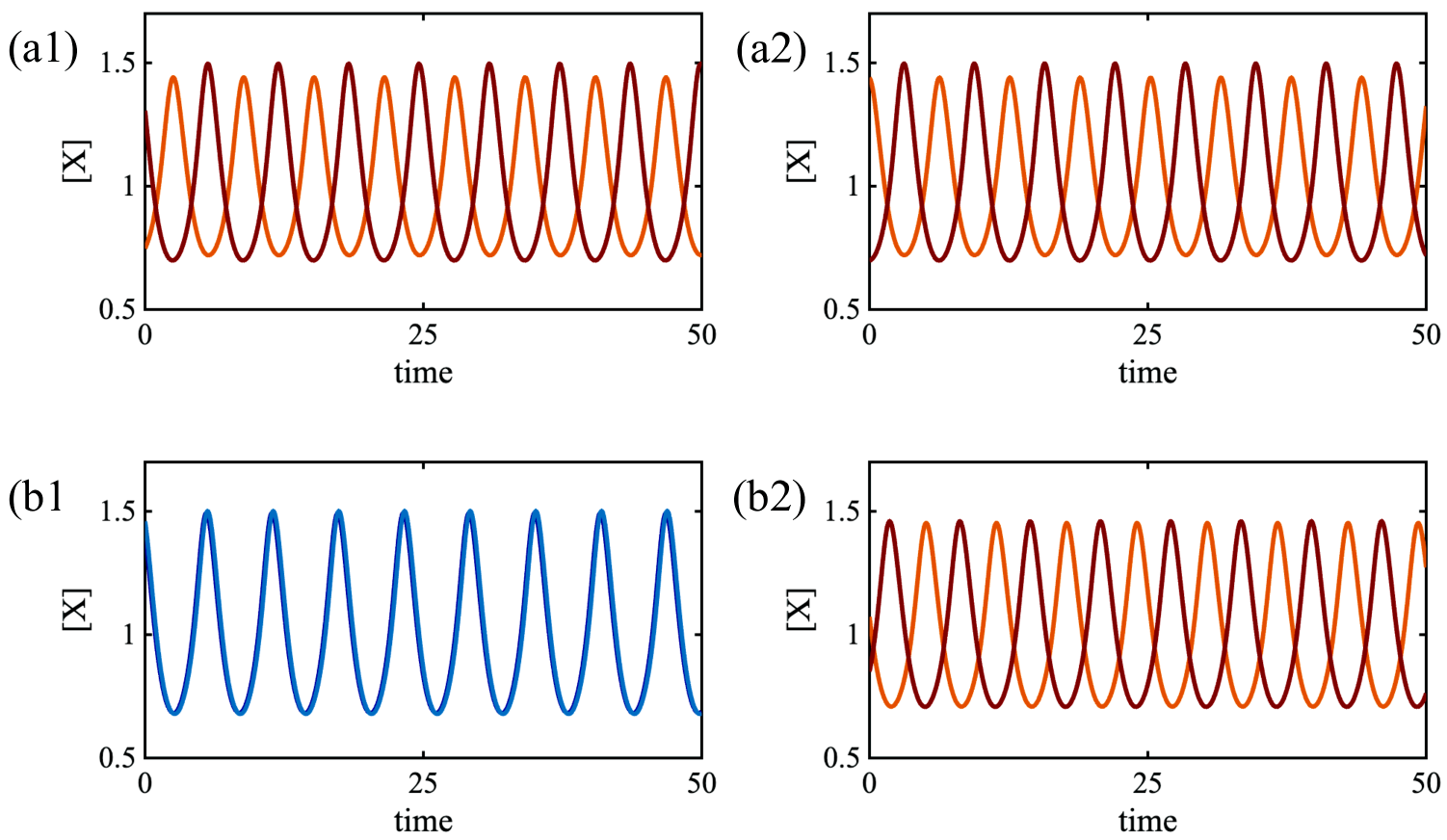}
    \caption{Bistability in delay-coupled Brusselator oscillators. Panels (a) and (b) show the behavior of weakly and strongly coupled systems, respectively. (a1-a2) $K = 0.01$, $\tau = 3.9$. Starting from both in-phase (a1) and anti-phase (a2) initial conditions, the system converges to an anti-phase synchronized state.  (b1-b2) $K = 0.05$, $\tau = 3.9$. Starting from in-phase initial conditions, the system remains in-phase synchronized (b1). Starting from anti-phase initial conditions, the system remains in anti-phase synchronized state (b2).}
    \label{fig:Supp1}
\end{figure}

Behavior of this system was studied at weak and strong coupling regimes, with a delay, $\tau = 3.9$, which is close to~$T/2$. The results are shown in Fig.~\ref{fig:Supp1}. In the weak coupling regime ($K = 0.01$), both in- and anti-phase initial conditions converge to anti-phase synchronized state. However, when coupling strength was increased to $K = 0.05$, in-phase initial conditions lead to in-phase synchrony and anti-phase initial conditions lead to anti-phase synchrony, exhibiting bistability. 

\section{Experimental Set-Up}
\label{app:Supp2}

\begin{figure}
    \centering
    \includegraphics[width=\linewidth]{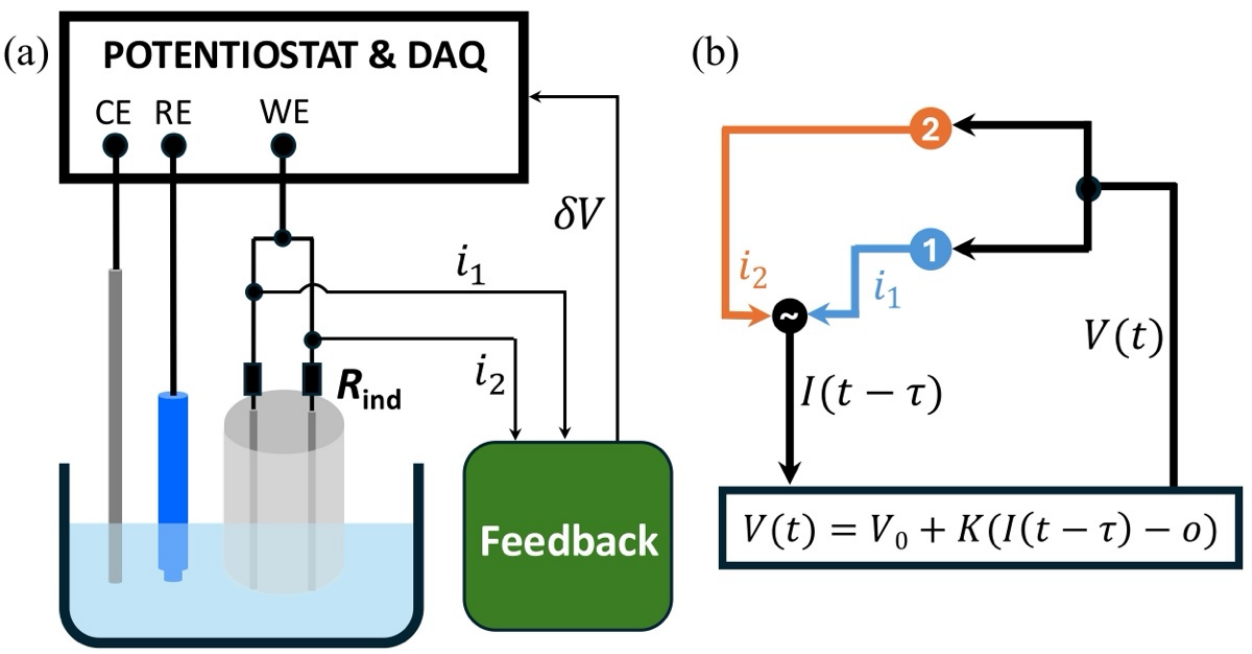}
    \caption{Figures/Experiments with nickel electrodissolution system. (a) Experimental setup. CE, RE, and WE are the counter, reference, and working electrodes. $R_\mathrm{ind}=1$ k\si{\ohm}, is the individual resistor. $i_1$ and $i_2$ are the the measured current of the wires at a given time~$t$ and~$\delta V$ is the perturbation applied to the circuit potential of the wires as a result of the feedback. (b) Feedback mechanism. $V(t)$ and $V_0(t)$ are the applied and offset circuit potentials, respectively, $K$ is the coupling strength, $I(t-\tau)=i_1(t-\tau)+i_2(t-\tau)$ is the  total current at a delay, $\tau$, and $o$ is an offset correction.}
    \label{fig:Supp2}
\end{figure}

The experimental setup (Fig.~\ref{fig:Supp2}(a)) consists of a three-electrode electrochemical cell with a Pt-coated Ti rod as counter (CE), a Hg/Hg\textsubscript{2}SO\textsubscript{4}/sat. K\textsubscript{2}SO\textsubscript{4} as reference (RE), and two Ni wires (Goodfellow Cambridge Ltd, 99.98\%, 1.0 mm diameter) embedded in epoxy as working electrodes (WE) connected to a potentiostat (ACM Instruments, Gill AC). The electrodes were immersed in a 3.0 M H\textsubscript{2}SO\textsubscript{4} solution as an electrolyte and kept at a constant temperature of 10\textsuperscript{o}C.

\begin{figure}
    \centering
    \includegraphics[width=1\linewidth]{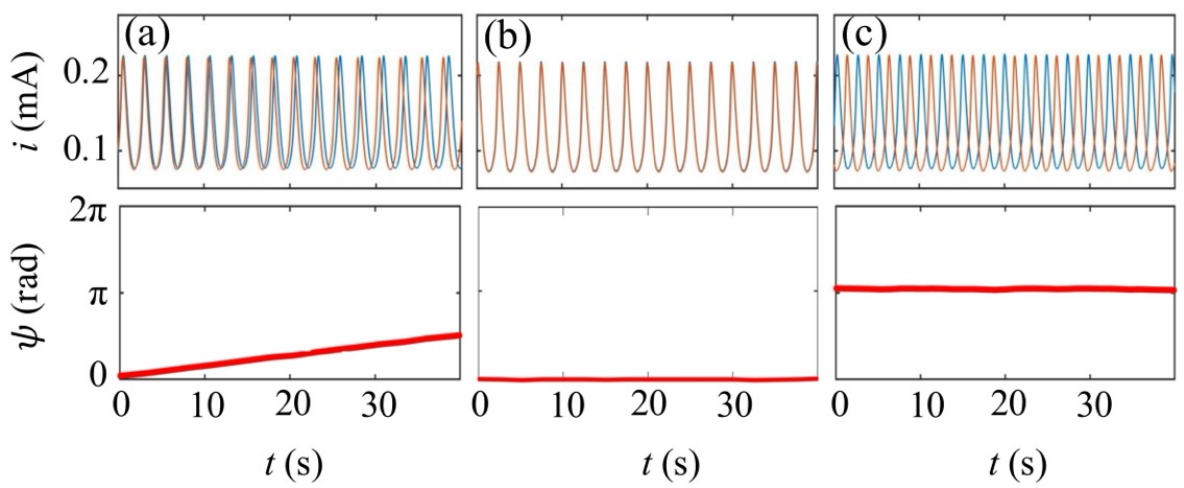}
    \caption{Oscillations at constant delay. (a) Uncoupled system ($K = 0$ VA$^{-1}$), current oscillations show no synchrony (top), the phase difference~$\psi$ extracted from the experimental data exhibits monotonic linear increase (bottom). (b) $\tau \approx 0$ ($K = -75$ VA$^{-1}$, $\tau = 0.00$ s), in-phase synchronized current oscillations (top), $\psi$~remains constant at $\psi=0$ (bottom). (c) $\tau \approx T/2$ ($K = -75$ VA$^{-1}$, $\tau = 0.92$ s), anti-phase synchronized current oscillations (top), $\psi$~remains constant at $\psi=\pi$ (bottom).}
    \label{fig:Supp3}
\end{figure}

A series resistance ($R_\mathrm{ind}=1$ k\si{\ohm}) was attached to each nickel wire. At the conditions mentioned above, the electrodissolution of nickel exhibits current oscillations due to the changing concentrations of species in the metal-electrolyte interface as the system goes through a Hopf bifurcation \cite{Haim1992ModelingPA}. In our experiments, Hopf bifurcation took place at circuit potential, $V_{0} = 1050$ mV. The experiments were were performed $30$ mV above this region at a constant circuit potential, $V_{0} = 1080$ mV, where the current shows periodic sinusoidal oscillations.

The two Ni wires were coupled with delayed linear feedback. The potentiostat interfaced with a real-time LabVIEW controller was used to measure the oscillating current ($i_k(t)$) of the two oscillators which were used to couple the elements following the coupling formula
\begin{equation}
\begin{aligned}
V(t) &= V_0 + K(I(t-\si{\tau}) - o)\\
\end{aligned}
\end{equation}
where~$V(t)$ and~$V_0(t)$ are the applied and offset circuit potentials, respectively, $K$~is the coupling strength, $I(t)=i_1(t)+i_2(t)$ is the  total current, $\tau$ is a time delay, and $o$~is an offset correction (set to the time-averaged total current). Fig.~\ref{fig:1}(b) shows a schematic illustration of the feedback mechanism.

From the currents of each wire,  the phases were extracted using peak-to-peak interpolation. When the electrodes are uncoupled, there is no synchrony between the two oscillators as observed in their current-time plots and the linearly drifting phase difference (Fig.~\ref{fig:Supp3}(a)). 
At $\tau = 0.00$ s, the oscillations are in-phase synchronized (see Fig.~\ref{fig:Supp3}(b)). As shown in Fig.~\ref{fig:Supp3}(c),
the two oscillators become anti-phase synchronized when the feedback delayed by about half of the average period of the oscillations ($\tau \approx T/2$).

\ 
\section*{Acknowledgments}

IZK and SC acknowledge support  by the National Institute of General Medical Sciences of the National Institutes of Health under award number R01GM157609, and by the National Institutes of Health under award number R01NS139415.

\section*{Code availability}

The computer algebra code is made available at \url{https://github.com/babettedewolff/reduction_SL}.

\bibliographystyle{apsrev4-2}
\def\urlprefix{}
\def\url#1{}
\bibliography{reduction}

\end{document}